\documentclass{article}

{}
{}
{}
\usepackage{color}
\usepackage{multirow}
\usepackage{makecell}
\usepackage{float}
\usepackage{setspace}
\usepackage{array}
\usepackage{longtable}
\usepackage{supertabular}
\usepackage{geometry}

\usepackage{mathtools}
\usepackage{amssymb}
\usepackage{graphicx}
\usepackage{lineno,hyperref}

\usepackage{amsfonts}
\usepackage{comment}

\begin{document}


\title{Locating All Real Solutions of Power Flow Equations: A Convex Optimization Based Method}

\author{Bin Liu\thanks{Department of Electrical Engineering, Tsinghua University, Beijing 100084, China and currently with School of Electrical Engineering and Telecommunications, The University of New South Wales, Sydney 2052, Australia. Email address: eeliubin@hotmail.com.}, Wei Wei\thanks{Department of Electrical Engineering, Tsinghua University, Beijing 100084, China. Email address: wei-wei04@mails.tsinghua.edu.cn.}, Feng Liu\thanks{Department of Electrical Engineering, Tsinghua University, Beijing 100084, China. Email address: lfeng@mail.tsinghua.edu.cn.}}

\maketitle

\begin{abstract}
This paper proposes a convex optimization based method that either locates all real roots of a set of power flow equations or declares no real solution exists in the given area. In the proposed method, solving the power flow equations is reformulated as a global optimization problem (GPF for short) that minimizes the sum of slack variables. All the global minima of GPF with a zero objective value have a one-to-one correspondence to the real roots of power flow equations. By solving a relaxed version of GPF over a hypercube, if the optimal value is strictly positive, there is no solution in this area and the hypercube is discarded. Otherwise the hypercube is further divided into smaller ones. This procedure repeats recursively until all the real roots are located in small enough hypercubes through the successive refinement of the feasible region embedded in a bisection paradigm. This method is desired in a number of power system security assessment applications, for instance, the transient stability analysis as well as voltage stability analysis, where the closest unstable equilibrium and all Type-I unstable equilibrium is required, respectively. The effectiveness of the proposed method is verified by analyzing several test systems. 
\end{abstract}

\section{Introduction}
Power flow (PF) equations formulate the relationships between complex bus voltages and nodal active/reactive power injections in the steady state. Solving PF equations is one of the most important and fundamental problems in power systems designing and operation both in theory and applications \cite{injp-1,cee-1,cee-2,injp-2,cee-3,injp-3}. PF solution also determines the initial condition of the power system when a transient behavior occurs. Due to the inherent non-linearity, PF equations can have many solutions. 

Most existing algorithms, such as the Newton-Raphson method \cite{Newton-PF,injp-4} and the fast-decoupled method \cite{FD-PF}, compute only one PF solution with the initial value around the nominal operating condition, i.e., the voltage magnitudes are near to 1 and the angles are near to 0. However, in a number of security assessment issues, more than a single solution are desired. For example, in power systems transient stability analysis, the closest unstable equilibrium point (Closest-UEP) is required to estimate the stability region \cite{CUEP-Method}. The Closest-UEP is closely related to the PF solution that has the lowest energy value of a so-called energy function among all others except the nominal operating point. Conceptually, to certify the Closest-UEP, one should evaluate and compare the values of the corresponding energy function at all these equilibriums. Several interesting methods have been proposed in \cite{CUEP-1,CUEP-2,CUEP-3} to compute the Closest-UEP without exploring all the solutions, but reported to have certain limitation in particular circumstances. Another example arises in voltage stability assessment \cite{Volt-Stability-1}, where multiple PF solutions are needed. Ref. \cite{Multi-PF-1} propose a method to compute a pair of PF solutions that are close to each other. Ref. \cite{Volt-Stability-2} reveals that the voltage instability phenomenon is closely related to the Type-1 equilibriums, which are defined as those where the Jacobian matrix of the PF equations has only one eigenvalue with positive real part. Ref. \cite{Type-1-UEP} proposes a continuation power flow based method to compute all Type-1 solutions. Ref. \cite{injp-5} presents an efficient starting process for calculating interval power flow solutions at maximum loading point while considering both load and line data uncertainties. 

To the best of our knowledge, finding all PF solutions is extremely challenging. Several research tries to address this problem, but few of them provides rigorous guarantee for the task. These methods include the continuation method \cite{All-Roots-Cont}, and the homotopy continuation method \cite{All-Roots-Homo-0,All-Roots-Homo-1,All-Roots-Homo-2}, both of which rely on tracing the solution paths with a parameter varying in a certain range. The first attempt was found in \cite{All-Roots-Homo-0}, where the homotopy continuation method in \cite{All-Roots-Homo-2} is adopted to find all PF solutions. However, the difficulty rests on its computational burden when the system size increases. Moreover, special care should be taken when the solution trajectory encounters a bifurcation. Ref. \cite{All-Roots-Cont} proposes a continuation algorithm to achieve this goal by exploring the topological structures of the solution set, while reducing the computation burden. However, a recent study \cite{All-Roots-Homo-1} shows the continuation algorithm in \cite{All-Roots-Cont} could miss some solutions in particular cases. A counter-example of a 5-bus system is given. Similar problem also exists in the method of \cite{Type-1-UEP}. Nevertheless, the homotopy continuation method in \cite{All-Roots-Homo-0} and \cite{All-Roots-Homo-3} is believed to be capable of finding all PF solutions \cite{All-Roots-Homo-1}, despite their poor scalability. A numerical polynomial homotopy continuation (NPHC) method is proposed in a latest publication \cite{All-Roots-Homo-2}. The NPHC method exploits the polynomial feature of PF equations, and is claimed to perform better than the traditional homotopy continuation method in terms of computational efficiency and ease of implementation, because no bifurcation is encountered. Instead, a large number of initial points should be supplied, whose number may grow exponentially with the problem size. In this regard, as all homotopy type methods, NPHC also suffers from the high computational burden. Nevertheless, as the solution trajectory starting at each initial point can be traced independently from all others, NPHC is easily parallelizable. The admissible problem size can go up to 13-bus system on a laptop, and 14-bus system on a a computing cluster with 64 processors. It should be pointed out that most of the computational effort of NPHC is spent in finding {\bf complex} solutions \cite{All-Roots-Homo-2}, which have no physical meaning or practical implication.

On the other hand, let's turn to another important problem in power system economic dispatch, the optimal power flow (OPF), which seeks the optimal generation portfolio subject to PF equations. There is no doubt that solving PF can be formulated as finding a global solution of a special OPF problem. However, due to the high non-convexity of the PF equations, this may not constitute an appealing approach if only a single solution is desired. During the past a few years, the convex relaxation based on semi-definite program (SDP) and second-order cone program (SOCP) technique have been proved to be very effective in computing the global solution of OPF problem \cite{CVX-OPF-1,CVX-OPF-2, CVX-OPF-3,CVX-OPF-4,CVX-OPF-5,CVX-OPF-6}, the relaxation is even tight under some mild conditions. 

Motivated by the recent advances in OPF and the method for solving general polynomial equations in \cite{Floudas}, this paper proposes a convex optimization based method to search all real PF solutions. We formulate the PF problem in rectangular coordinate as a quadratic constrained programming (QCP), whose objective is to minimize the sum of slack variables. Each PF solution corresponds to a global minimum with a zero objective value. In \cite{Floudas}, only the convex-concave envelop is adopted to relax the corresponding polynomial optimization problem into a linear program. In this paper, we combine the convex relaxation (including both SDP and SOCP) approach with the convex-concave envelops, yielding an enhanced convex relaxation of the QCP associated with PF equations. We use the bisection search scheme in \cite{Floudas} to locate all real PF roots in separated small enough rectangles by successively narrowing the upper and lower bounds of the decision variables, thus the relaxation is gradually tightened. We will demonstrate that the convex constraint plays an important role on reducing the computation time when the system size grows larger. It should also be mentioned that we do not require the relaxation always be tight, which is however, desired by OPF applications. In contrast to all the continuation inspired methods, our method relies on commercial SDP or SOCP solvers, so is easy to implement from a computational point of view. Moreover, it can reliably locate all the {\bf real} solutions or declare infeasibility of PF problems in given areas.

The rests of this paper are organized as follows. In Section 2, we review the PF equations in rectangular form and formulate it as a QCP. In section III, we present the convex optimization based bisection search algorithm. In section IV, we provide the computation results on several test systems, and compare the performance of the algorithm with and without convex constraint. Experiments on a 9-bus and 14-bus system illustrates its flexibility in detecting infeasibility.

\section{Mathematic Formulation}
Symbols and notations used in this section will be defined first. A power network is modeled by $\mathcal{N}=(\mathcal{B,L})$, where $\mathcal{B}=\mathcal{PQ}\cup \mathcal{PV} \cup \mathcal{V \theta} $ is the set of buses with the number of buses as $|\mathcal{B}|=n$, $\mathcal{PQ}$ denotes the set of buses whose active/reactive power injection is constant, $\mathcal{PV}$ denotes the set of buses whose active power injection and voltage magnitude is constant, $\mathcal{V\theta}$ denotes the reference bus whose complex voltage is given and $\mathcal{L}$ is the set of transmission lines. Let the nodal admittance matrix be $Y\in {\mathbb C}^{n\times n}$, where ${\mathbb C}^{n\times n}$ is the set of $n \times n$ complex matrix, and its elements be $Y_{ij} =G_{ij} +\mbox{j}B_{ij}(\forall (ij)\in \mathcal{L})$, where $\mbox{j}=\sqrt {-1}$. Let $p_i^g /q_i^g $ be the active/reactive power of generator at bus $i$ ($p_i^g =q_i^g =0$ if there is no generator at that bus), and $p_i^d /q_i^d $ be the active/reactive power demand at bus $i$ ($p_i^d =q_i^d =0$ if there is no load at that bus). Denoting the voltage magnitude and phase angle of bus $i$ as $|V_i|$ and $\theta_i$, respectively, its complex voltage in the rectangular coordinate can be expressed as $V_i = e_i +\mbox{j}f_i=|V_i|\cos{\theta_i}+\mbox{j}|V_i|\sin{\theta_i}$. With the above notations, the PF equations are given in the rectangular form as follows.
\begin{eqnarray}\label{eq-1-1}
p_i^g - p_i^d = \sum\limits_{j=1}^n {[G_{ij} (e_i e_j  + f_i f_j ) - B_{ij} (e_i f_j -e_j f_i )]}\\
\quad \forall i \in \mathcal{PV \cup PQ}\nonumber\\      
\label{eq-1-2}
q_i^g - q_i^d = \sum\limits_{j=1}^n {[-B_{ij} (e_i e_j + f_i f_j ) - G_{ij} (e_i f_j -e_j f_i )]}\\
\quad \forall i\in \mathcal{PQ}\nonumber
\end{eqnarray}
\begin{eqnarray}
\label{eq-1-3}e_i^2 +f_i^2 =| {V_i } |^2,\ \forall i \in \mathcal{PV}\\
\label{eq-1-4}e_i = |V_i|, \ f_i = 0,\forall i \in \mathcal{V\theta}
\end{eqnarray}

Equation set \eqref{eq-1-1}-\eqref{eq-1-4} contains $2n$ equations, and $2n$ variables as well. For notation brevity, introduce the standard basis vectors in ${\mathbb R}^n$ as $b_1$, $b_2$, $\cdots$, $b_n$, for each basis vector  $b_k(k\in \{1,2,\cdots ,n\}$, the $k$-th element of $b_k$ is 1, and the rest ones are 0. Define the following matrices
\begin{eqnarray}
Y_k =b_k b_k^T Y \in {\mathbb C}^{n\times n}\nonumber\\
M_k =\left[\begin{aligned}
b_kb_k^T&~&0 \\
0&~&b_k b_k^T
\end{aligned}\right]\in {\mathbb R}^{2n\times 2n}\nonumber\\
Z_k =\frac{1}{2}\left[\begin{aligned}
{\mbox{Re}\{Y_k +Y_k^T \}}  &~& {\mbox{Im}\{Y_k^T -Y_k \}}  \\
{\mbox{Im}\{Y_k -Y_k^T \}}  &~& {\mbox{Re}\{Y_k +Y_k^T \}}  \\
\end{aligned}\right]\in {\mathbb R}^{2n\times 2n}\nonumber\\
\bar Z _k =-\frac{1}{2} \left[\begin{aligned}
{\mbox{Im}\{Y_k +Y_k^T \}}  &~& {\mbox{Re}\{Y_k -Y_k^T \}}  \\
{\mbox{Re}\{Y_k^T -Y_k \}}  &~& {\mbox{Im}\{Y_k +Y_k^T \}}  \\
\end{aligned}\right]\in {\mathbb R}^{2n\times 2n}\nonumber
\end{eqnarray}
where operator $\mbox{Re}\{\cdot\}/\mbox{Im}\{\cdot\}$ represents the real/imaginary part of a complex vector or matrix. 

The PF equations \eqref{eq-1-1}-\eqref{eq-1-4} can be arranged in a compact form as follows \cite{CVX-OPF-1} using above notations(referred as \textbf{PF})
\begin{eqnarray}
\label{eq-2-1} x^TZ_k x=p^{ \text{in}}_k,~~\forall k \in \mathcal{PV \cup PQ}  \\ 
\label{eq-2-2}x^T\bar Z_k x = q^{ \text{in}}_k,~~\forall k \in \mathcal{PQ}   \\ 
\label{eq-2-3}x^TM_k x=V_k^2, ~~ \forall k \in \mathcal{PV}      \\
\label{eq-2-4}x_k = |V_k|, x_{k+n} = 0,~~k \in \mathcal{V\theta}          
\end{eqnarray}
where $x=[e^T \ f^T]^T$ is the variable, $p^{ \text{in}}_k = p_k^g -p_k^d(\forall k \in \mathcal{PV \cup PQ})$ and $q^{ \text{in}}_k =q_k^g -q_k^d(\forall k \in \mathcal{PQ})$ are active and reactive power injections, respectively.  

Next we analyze the range of variables. The voltage magnitudes $|V_k|(\forall k \in \mathcal{PV})$ are given by power flow data. From equation \eqref{eq-1-3}, it is clear that $e_k$ and $f_k$ of $PV$ buses, i.e. buses belonging to $\mathcal{PV}$, are bounded by $|V_k|$. As for the $PQ$ buses, i.e. buses belonging to $\mathcal{PQ}$, their voltage magnitudes are also bounded, but the value $|V_k|=\sqrt {e^2_k+f^2_k}$ is not instantly clear at hand. The predictor-corrector technique proposed in \cite{PF-Boundary} can be applied to find such boundaries if desired. In the paper, we simply choose $V_k^{\text{min}}=0,V_k^{\text{max}}=1.5$ as the lower and upper values of magnitude for $PQ$ buses, $\theta_k^{\text{max}}=-\theta_k^{\text{min}}=\pi$ as the lower and upper values of phase angle both $PV$ and $PQ$ buses. Note that this is only a boundary estimation for the voltage magnitude, we do not actually impose voltage magnitude constraints on $PQ$ buses. In summary, we have the following bound constraints for variable $x$ 
\begin{eqnarray}
\label{eq-3}[x^l_k;~x^l_{n+k}]\le [x_k;~x_{n+k}] \le [x^u_k;~x^u_{n+k}]
\end{eqnarray}
where
\begin{eqnarray}
\label{eq-4-1}x_k^u=x_{n+k}^u=-x_k^l=-x_{n+k}^l=V_k^{\text{max}},~\forall k\in \mathcal{PQ}\\
\label{eq-4-2}x_k^u=x_{n+k}^u=-x_k^l=-x_{n+k}^l=|V_k|,~\forall k\in \mathcal{PV}\\
\label{eq-4-3}x_k^u=x_k^l=|V_k|,x_{n+k}^u=x_{n+k}^l=0,~\forall k\in \mathcal{V\theta}
\end{eqnarray}

By introducing slack variables $s_{pk}^+$, $s_{pk}^-$, $s_{qk}^+$, $s_{qk}^-$, $s_{vk}^+$, $s_{vk}^-$, PF equations \eqref{eq-2-1}-\eqref{eq-2-4} can be written as the following optimization problem, referred as (\textbf{GPF}) in the following context
\begin{eqnarray}
\label{eq-5-1}S_{opt} =\min \sum\limits_{k\in \mathcal{PV \cup PQ}} {(s_{pk}^+ +s_{pk}^- )}
+\sum\limits_{k\in \mathcal{PQ}} {(s_{qk}^+ +s_{qk}^- )} \\
+\sum\limits_{k\in \mathcal{PV}}
{(s_{vk}^+ +s_{vk}^- )}  \\
\label{eq-5-2}\mbox{tr}(XZ_k )+s_{pk}^+ -s_{pk}^- =p^{\text{in}}_k,\ \forall k \in
\mathcal{PV \cup PQ}    \\
\label{eq-5-3}\mbox{tr}(X\bar Z _k )+s_{qk}^+ -s_{qk}^- =q^{\text{in}}_k,\ \forall k \in 
\mathcal{PQ}   \\
\label{eq-5-4}\mbox{tr}(XM_k )+s_{vk}^+ -s_{vk}^- =|V_k|^2,\ \forall k \in \mathcal{PV}  \\
\label{eq-5-5}s_{pk}^+ ,s_{pk}^- \ge 0,\  \forall k\in \mathcal{PV \cup PQ}  \\
\label{eq-5-6}s_{qk}^+ ,s_{qk}^- \ge 0,\  \forall k\in \mathcal{PQ}  \\
\label{eq-5-7}s_{vk}^+ ,s_{vk}^- \ge 0,\  \forall k\in \mathcal{PV}  \\
\label{eq-5-8}X = x x^T                                    \\
\label{eq-5-9}x \in {\mbox {BOX}(x^l,x^u)}      
\end{eqnarray}
where ${\mbox {tr}} (\cdot)$ represents the matrix trace operator, the hypercube ${\mbox {BOX}(x^l,x^u)} = \{x|x^l\le x \le x^u\}$ with the bound parameter $x^l$ and $x^u$.

The relationship between {\bf PF} and {\bf GPF} is summarized as follows: each global minima $(X^\ast,x^\ast,s^\ast )$ with a zero optimal value of {\bf GPF} corresponds to a solution of {\bf PF}, and vice versa. Moreover, if the global minimum $S_{opt}>0$ then {\bf PF} has no solutions in the rectangle BOX($x^l, x^u$). Therefore, locating all power flow solutions requires finding all the global optimal solutions of {\bf GPF} with a zero-objective value, which is very challenging because {\bf GPF} is non-covnex due to constraint \eqref{eq-5-8}. Any local optimization method is not able to find all global minima of {\bf GPF} with rigorous guarantee or even deduce a wrong statement that {\bf PF} has no solution. 

In some circumstances, the solutions of power flow problem falling into a more restricted area may be interested in and the corresponding constraints can be easily incorporated. One of them is the phase angle difference (PAD) constraint for each line, which can be formulated as the following constraints.
\begin{equation}\label{eq-6}\begin{split}
-(X_{i,j}+X_{n+i,n+j})\tan{\Delta\theta_{ij}^{\text{max}}}\le X_{j,n+i}-X_{i,n+j}\\
\le (X_{i,j}+X_{n+i,n+j})\tan{\Delta\theta_{ij}^{\text{max}}},~~\forall (ij)\in \mathcal{L}
\end{split}\end{equation} 
where $\Delta\theta_{ij}, \Delta\theta_{ij}^{\text{max}}$ are the PAD of line $(ij)$ and its upper bound, respectively.

Another one is related to when more specified values of $V_k^{\text{min}}$, $V_k^{\text{max}}$, $\theta_k^{\text{max}}$, $\theta_k^{\text{min}}$ are given. The corresponding constraints will depend on the given values of the parameters. For example, the constraints of a $PQ$ bus with $0\le V_k^{\text{min}}\le V_k^{\text{max}}$ and $0\le\theta_k^{\text{min}}\le\theta_k^{\text{max}}\le\pi/2$ are showed as \eqref{eq-7}, which will replace \eqref{eq-4-1} in the formulated models.
\begin{equation}\label{eq-7}\begin{split}
x^l_k=V_k^{\text{min}}\cos{\theta_k^{\text{max}}},x^u_k=V_k^{\text{max}}\cos{\theta_k^{\text{min}}}\\
x^l_{n+k}=V_k^{\text{min}}\sin{\theta_k^{\text{min}}},x^u_{n+k}=V_k^{\text{max}}\sin{\theta_k^{\text{max}}}\\
\end{split}\end{equation} 

The constrains can be formulated similarly for $PV$ buses when the boundary values are specified, which are omitted here. However, if the boundary values of any bus are not specified, we still estimate that  $V_k^{\text{min}}=0,V_k^{\text{max}}=1.5$ for $PQ$ buses and $\theta_k^{\text{max}}=-\theta_k^{\text{min}}=\pi$ for both $PQ$ and $PV$ buses.

\section{Methodology}
Based on the methodology proposed in \cite{Floudas}, SDP relaxation proposed in \cite{CVX-OPF-1} and SOCP relaxation proposed in \cite{CVX-OPF-2}, we outline an enhanced algorithm to locate all the global solutions of ({\bf{GPF}}) through convex relaxation in box-constrained area and successively tightening the bounds of hypercubes that may contain a solution of {\bf PF}. 

\subsection{Convex Relaxation of {\bf GPF}}

By noticing the fact that for any vector $v$ with proper dimension, $v^T Xv = (v^T x)(x^T v) = (v^T x)^2\ge 0$ means matrix variable $X$ should be positive semi-definite. In this regard, the standard SDP relaxation of {\bf GPF} is to replace constraint \eqref{eq-5-8} with the following matrix inequality \cite{CVX-OPF-1}
\begin{equation}
\label{eq-8}X \succeq 0
\end{equation}
where $\succeq$ represents matrix $X$ is positive semi-definite. This paradigm is widely adopted in relaxing the traditional optimal power flow problem, and is receiving increasing attention nowadays. Constraint \eqref{eq-8} is only a necessary consequence of constraint \eqref{eq-5-8}. Sufficiency holds if $\mbox {rank} \{X\} =1$. However, rank constraint is non-convex. To further tighten the relaxation and acquire stronger bounds, we will utilize valid linear inequalities, i.e., the reformulation-linearization technique (RLT) developed in \cite{RLT}. RLT is motivated from the following fact: $\forall i,j \in \{ 1,2, \cdots n\} ,i \le j$ and $\forall x \in {\mbox {BOX}}(x^l,x^u)$, the following equalities hold
\begin{equation}\label{eq-9}\begin{split} 
(x_i - x^l_i)(x_j - x^l_j) = X_{i,j} - x_i x^l_j - x_j x^l_i + x^l_i x^l_j \\
(x_i - x^l_i)(x^u_j - x_j) = x_i x^u_j - X_{i,j} + x_j x^l_i - x^l_i x^u_j  \\
(x^u_i - x_i)(x_j - x^l_j) = x^u_i x_j - x^u_i x^l_j - X_{i,j} + x_i x^l_j \\
(x^u_i - x_i)(x^u_j - x_j) = x^u_i x^u_j - x^u_i x_j - x_i x^u_j + X_{i,j} 
\end{split}
\end{equation}

Equation \eqref{eq-9} reveals additional connections between $x$ and $X$ indicated by the lower/upper bound $x^l/x^u$. In this section, we will adopt its compact form in \cite{RLT-Comp}. Define RLT constraints for variable $(x,X)$ with parameter $x^l$ and $x^u$ as follows
\begin{equation}\label{eq-10}\begin{split}
&\mbox {RLT}(x^l,x^u)\\
&=\left\{\makecell{(x,X)	\in \\{\mbox {BOX}(x^l,x^u)} \times \mathbb S^n} \left| 
\begin{array}{c}
{x^l x^T + x (x^l)^T - x^l (x^l)^T \le X}  \\
{x^u x^T + x (x^u)^T - x^u (x^u)^T \le X}  \\
{X \le x (x^u)^T + x^l x^T- x^l (x^u)^T}  \\
\end{array} \right. \right\}
\end{split}\end{equation}
where ${\mathbb S}^n$ is the set of real symmetric matrix.

The convex relaxation of {\bf GPF} combining both SDP constraint \eqref{eq-8} and RLT constraint \eqref{eq-10} in the rectangle BOX($x^l, x^u$) is
\begin{equation}\label{eq-11}
S_{sdp} =\min{\{\eqref{eq-5-1}|s.t. (x,X) \in \mbox {RLT}(x^l,x^u),\eqref{eq-5-2}-\eqref{eq-5-7},\eqref{eq-8}\}}
\end{equation}
Problem \eqref{eq-11} is refered as ${\bf GPF_{sdp}}$ in the following context. As all constraints in \eqref{eq-5-2}-\eqref{eq-5-7} and RLT($x^l,x^u$) are linear, ${\bf GPF_{sdp}}$ is a SDP problem, whose global optimal solution can be efficiently computed by commercial solvers. Compared with the standard SDP relaxation in \cite{CVX-OPF-1}, an important feature of ${\bf GPF_{sdp}}$ is that the gap between the optimal values of ${\bf GPF_{sdp}}$ and {\bf GPF} is becoming smaller when the upper and lower bound is getting closer \cite{Floudas}. Compared with the pure RLT relaxation, the SDP constraint strengthen the relaxation in each rectangle, and can reduce the computation time remarkably. 

The SDP constraint \eqref{eq-8} can be further relaxed to a set of SOCP constraints \cite{CVX-OPF-2}
\begin{equation}\label{eq-12}\begin{split}
(X_{i,i}+X_{n+i,n+j})(X_{j,j}+X_{n+j,n+j})\le(X_{i,j}+X_{n+i,n+j})^2\\
+(X_{n+i,j}+X_{i,n+j})^2,~\forall (ij)\in \mathcal{L}
\end{split}\end{equation}
which leads to another convex relaxation form of {\bf GPF} given below. 
\begin{equation}\label{eq-13}
S_{socp} =\min{\{\eqref{eq-5-1}|s.t. (x,X) \in \mbox {RLT}(x^l,x^u),\eqref{eq-5-2}-\eqref{eq-5-7},\eqref{eq-12}\}}
\end{equation}

Problem \eqref{eq-13}, which is a SOCP problem, is refereed as ${\bf GPF_{socp}}$ in the following text. It is noteworthy that for power flow problem of radial power networks, ${\bf GPF_{socp}}$ is equivalent to ${\bf GPF_{sdp}}$ according to \cite{CVX-OPF-2}. 

Besides, the convex relaxed {\bf GPF} problem with only RLT constraints is refereed as ${\bf GPF_{rlt}}$ and is formulated as follows for the convenience of case studies.  
\begin{equation}\label{eq-14}
S_{rlt} =\min{\{\eqref{eq-5-1}|s.t. (x,X) \in \mbox {RLT}(x^l,x^u),\eqref{eq-5-2}-\eqref{eq-5-7}\}}
\end{equation}

According to the convex relaxation methodology employed, ${\bf GPF_{sdp}}$ is the most tight form of {\bf GPF} while ${\bf GPF_{rlt}}$ has the least tightness, which leads to the following inequalities.
\begin{equation}\label{eq-7+}
S_{opt}\ge S_{sdp}\ge S_{socp} \ge S_{rlt}
\end{equation}

\subsection{The algorithm}
The theoretical basis of the algorithm is explained below. Provided the region BOX($x^l,x^u$) under investigated, solve the convex relaxed {\bf GPF} (${\bf GPF_{cvx}}$ for short) problem, i.e. ${\bf GPF_{rlt}}$, ${\bf GPF_{socp}}$ or ${\bf GPF_{sdp}}$. For simplicity, $S_{cvx}$ is used to denote the optimal value of the convex relaxed problem, which equals to $S_{rlt}$, $S_{socp}$ and $S_{sdp}$ for ${\bf GPF_{rlt}}$, ${\bf GPF_{socp}}$ and ${\bf GPF_{sdp}}$ respectively. If the optimal value $S_{cvx>0}$, we can conclude that {\bf PF} has no solution because $S_{cvx}$ is a lower bound of the optimal value of {\bf GPF} according to \eqref{eq-7+}, or else if the optimal value $S_{cvx}= 0$, there may be none, one, or multiple solutions in BOX($x^l,x^u$), and we further divide the current hypercube into smaller sub-hypercubes and solve ${\bf GPF_{cvx}}$ in each of them. This procedure is repeated recursively, until all global solutions are located in small enough rectangles (thus can be refined by standard PF algorithm), and the remaining area is certified to contain no solution. The algorithm is formally provided below.

{\bf Algorithm 1}

\textbf{Step 1}: Initialize error tolerance $\varepsilon_R$, $\varepsilon _V$, the set of candidate hypercubes $B_C = \emptyset$, the set of solutions ${\rm O}^\ast =\emptyset$, and the current set of bounds $L=\{x^l\}$, $U=\{x^u\}$.

\textbf{Step 2}: Set $l=L_1,u=U_1$ (where $L_1/U_1$ represents the first element of set $L/U)$, solve ${\bf GPF_{cvx}}$ in BOX($l,u$), the optimal solution is $x^\ast$ and the optimal value is $S_{cvx}$.

1) If $S_{cvx}>\varepsilon_R$, update $U=U\backslash U_1,L=L\backslash L_1$.

2) If $S_{cvx}\le \varepsilon_R$ and $\left\| u - l \right\|_\infty \le \varepsilon_V$, BOX($l,u$) contains one solution, update $B_C = B_C \cup \mbox {BOX}(l,u)$, ${\rm O}^\ast ={\rm O}^\ast \cup \{x^\ast \}$, update  $U=U\backslash U_1,L=L\backslash L_1$.

3) If $S_{cvx}\le \varepsilon_R$ and $\left\| u - l \right\|_\infty 
>\varepsilon_V$, update $U=U\backslash U_1,L=L\backslash L_1$, find one entry $k$ such that $u_k -l_k =\left\| {u-l} \right\|_\infty $, the current rectangle is partitioned into two smaller ones indicated by the following two sets of boundaries
\[l^1 = l, \quad  u^1 = [u_1,\cdots,(l_k + u_k)/2,\cdots,u_{2n}]^T \]
\[u^2 = u, \quad  l^2 = [l_1,\cdots,(l_k + u_k)/2,\cdots,l_{2n}]^T\]
update $L=L\cup l^1 \cup l^2$ and $U = U \cup u^1 \cup u^2$ accordingly. 

\textbf{Step 3}: If $U \ne \emptyset$, go to {\textbf{Step 2}}.

\textbf{Step 4}: Retrieve the sets $B_C$ and $\rm {O}^\ast$, for each hypercube in $B_C$, solving \textbf{PF} using the traditional power flow algorithm with the corresponding initial point in ${\rm O}^\ast $ , and 
record all solutions.

Some additional discussions are provided.

1) $\varepsilon _V $ is the tolerance to judge whether $l$ and $u$ is close enough and has a significant impact on the computational efficiency and accuracy of Algorithm~1. There are two considerations when selecting this parameter: a) In general, $\varepsilon _V $ should be less than the distance between the closest two solutions, such that each rectangle in $B_C$ contains at most one solution; b) from an accuracy point of view, $\varepsilon _V $ should be as small as possible such that the relaxation is tight enough. However, on the one hand, one does not have the information on the distribution of power flow solutions in advance; on the other hand, the algorithm will suffer from computational burden due to the extensive branching process in Step 2 when $\varepsilon _V $ becomes too small. If we can estimate the distance $d_P $ between the closest two solutions (to the best of our knowledge, this is still an open problem), we can choose $\varepsilon _V =d_P $ and call Algorithm 1 to find all the rectangles that contains only one solution, and further apply Algorithm 1 in each area using a smaller $\varepsilon_V$ to refine the result, or using other methods to find the only one solution. Nevertheless, we can also specify a less theoretical $\varepsilon _V$ from experience in the beginning, and reduce its value in further computation process. In this regard, Algorithm~1 should be called recursively.   This trick is important when the system size grows larger. It identifies the area without a solution quickly without further branching. Here we assume when the rectangles in $B_C$ are small, the initial values provided by ${\rm O}^\ast $ are good enough such that the traditional power flow algorithm can  converge.  

2) In \cite{Floudas}, only the RLT constraint (5) is adopted to obtain a linear programming relaxation of polynomial equations. In this paper, the combination of RLT and convex constraints can tighten the relaxation so as to enhance the efficiency. The convex constraint is valuable especially when the problem size grows larger, although solving convex problem is more expensive than solving linear programming problem. This will also be demonstrated in case studies. Beside the convex constraints, other relaxation technique is also valid for ${\bf GPF_{cvx}}$ , such as other valid linear inequalities and second order cone inequalities in \cite{RLT-Comp}. From a computational point of view, solving SOCPs is less expensive than solving SDPs, on the other hand, SOCP relaxation of OPF is usually weaker than the SDP relaxation for mesh networks \cite{Strong-SOCP}. It is noteworthy that we do not require the relaxation always be tight, which is different from OPF applications, because the gap between the optimal values of {\bf GPF} and ${\bf GPF_{cvx}}$  would be small enough with the bound $l$ and $u$ getting closer.

3) Our algorithm would be more efficient in certifying that a set of power flow equations have no real solution, because the optimal value of {\bf GPF} is strict positive over the rectangle BOX($x^l,x^u$) in such circumstance, thus Algorithm 1 is expected to encounter fewer branching actions in Step 2. So far as we know, the homotopy continuation method has to compute all complex solution, so as to draw such a conclusion. Our algorithm is also easily parallelizable by partitioning the initial rectangle into disjoint smaller subregions.

\section{Case Studies}
In this section, we apply the proposed method on the 5-bus system and 7-bus system investigated in \cite{All-Roots-Homo-0,All-Roots-Homo-1} for the purpose of validating the method by comparing its outcome with known results. Further studies on IEEE 9-bus and 14-bus systems show its ability to detect infeasibility or special PF solution under various settings. All the tests are implemented on a laptop computer with Intel i7-3520M CPU and 8 GB memory. LPs are solved by CPLEX 12.5. SDPs and SOCPs are solved by MOSEK.

\subsection{5-Bus and 7-Bus Systems}
The power flow solutions of the 5-bus and 7-bus systems reported by the proposed method are shown through Table \ref{tab-1} and \ref{tab-2}. According to the results in \cite{All-Roots-Homo-0,All-Roots-Homo-1,All-Roots-Homo-2}, all the solutions have been found. Moreover, we test the computational efficiency by varying $\varepsilon _V $ from $10^{-4}$ to $10^{-1}$ with $\varepsilon _R = 10^{-5}$. The computation time is shown in Table \ref{tab-3}. 
\begin{table}[ht]\renewcommand\arraystretch{1.0}
	\caption{Solutions of the 5-Bus System}
	\vspace{6pt}
	\centering
	\begin{tabular}{>{\scriptsize}c|>{\scriptsize}c|>{\scriptsize}c|>{\scriptsize}c|>{\scriptsize}c|>{\scriptsize}c}
		\hline
		Solution & $\delta_1$ & $\delta_2$ & $\delta_3$ & $\delta_4$ & $\delta_5$\\
		\hline
		1 &    0   &   1.286    &   22.06    &    2.194   &   0.372  \\
		2 &    0   &   2.187    &   45.92   &    46.62   &  -144.00  \\
		3 &    0   &   0.166    &   171.20    &    0.028   &  -0.710  \\
		4 &    0   &  -0.897    &  -168.40    &    44.39   &  -145.10  \\
		5 &    0   &  -169.40    &  -10.99    &    165.90   &  -25.38  \\
		6 &    0   &  -168.70    &   3.182    &   -167.10   &  -167.90  \\
		7 &    0   &  -171.40    &  -99.23    &    50.72   &  -141.80  \\
		8 &    0   &  -168.70    &  -172.90    &    44.01   &  -145.30  \\
		9 &    0   &  -169.30    &  -160.90    &    166.10   &  -22.90  \\
		10 &   0   &  -169.90    &  -148.20    &   -167.10   &  -168.90  \\
		\hline
	\end{tabular}  
	\label{tab-1}
\end{table}
\begin{table}[ht]\renewcommand\arraystretch{1.0}
	\caption{Solutions of the 7-Bus System}
	\vspace{6pt}
	\centering
	\begin{tabular}{>{\scriptsize}c|>{\scriptsize}c|>{\scriptsize}c|>{\scriptsize}c|>{\scriptsize}c}
		\hline
		solution  &      1     &      2     &       3    &       4    \\
		\hline
		$V_1$   &      1.000     &      1.000     &       1.000    &          1.000 \\
		$V_2$   &      1.089 &      0.718 &      0.354 &      0.281 \\
		$V_3$   &      0.977 &      0.573 &      0.419 &       0.560 \\
		$V_4$   &      0.934 &      0.166 &      0.223 &      0.595 \\
		$V_5$   &      0.947 &      0.394 &      0.406 &      0.679 \\
		$V_6$   &      0.973 &      0.713 &      0.675 &      0.789 \\
		$V_7$   &      0.979 &      0.652 &      0.537 &      0.655 \\
		\hline
		$\delta_1$ &          0 &          0 &          0 &          0 \\
		$\delta_2$ &      4.918 &     15.406 &      84.980 &    102.957 \\
		$\delta_3$ &      -3.090 &     -5.719 &     -7.358 &     -6.538 \\
		$\delta_4$ &     -8.644 &    -62.623 &    -55.067 &    -19.191 \\
		$\delta_5$ &     -5.925 &    -15.722 &    -18.189 &     -11.290 \\
		$\delta_6$ &     -2.561 &     -3.431 &     -4.139 &     -4.034 \\
		$\delta_7$ &     -2.716 &     -4.651 &     -5.641 &       -5.200 \\
		\hline
	\end{tabular}  
	\label{tab-2}
\end{table}

\begin{figure}[ht]
	\centering
	\includegraphics[scale=0.25]{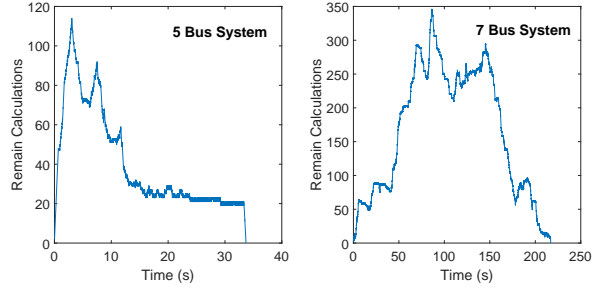}
	\caption{Number of hypercubes needed exploration during computation of ${\bf GPF_{sdp}}$ with $\varepsilon_V = 10^{-2}$}
	\label{fig-1}
\end{figure}

Overall, the computation time decreases with $\varepsilon _V $ increasing. For the 5-bus system, the ${\bf GPF_{rlt}}$ relaxation is more efficient. For the 7-bus system, both ${\bf GPF_{socp}}$ and ${\bf GPF_{sdp}}$ appear to be significantly more efficient than ${\bf GPF_{rlt}}$, demonstrating that convex constraint \eqref{eq-8} or \eqref{eq-12} plays an important role in tightening the relaxation and accelerating convergence. In Algorithm 1, the initial hypercube is divided into smaller hypercubes that may contain a solution. The number of hypercubes need exploration during computation of ${\bf GPF_{sdp}}$ is illustrated in Figure \ref{fig-1}.
\begin{table}[ht]\renewcommand\arraystretch{1.0}
	\caption{Computation time comparison}
	\vspace{6pt}
	\centering
	\begin{tabular}{>{\scriptsize}c|>{\scriptsize}c|>{\scriptsize}c|>{\scriptsize}c|>{\scriptsize}c}
		\hline
		\multicolumn{2}{c|}{{\scriptsize Computation Time (s)}}    &    $\varepsilon_V=10^{-4}$   &  
		$\varepsilon_V=10^{-2}$   &    $\varepsilon_V=10^{-1}$   \\
		\hline
		5-bus  &  ${\bf GPF_{rlt}}$ &    37.92  &     22.51 &    16.31  \\
		system & ${\bf GPF_{socp}}$  &    64.72  &     31.60 &    22.37  \\
		& ${\bf GPF_{sdp}}$  &    115.59  &     33.70 &    24.44  \\
		\hline
		7-bus &  ${\bf GPF_{rlt}}$   &    6519.72   &     6349.52  &    6167.67   \\
		system & ${\bf GPF_{socp}}$ &    449.86  &     367.54 &    341.02  \\
		& ${\bf GPF_{sdp}}$ &    367.32  &    209.04 &    174.82  \\
		\hline
	\end{tabular} 
	\label{tab-3} 
\end{table}

\subsection{9-Bus System}
\begin{figure}
	\centering
	\includegraphics[scale=0.25]{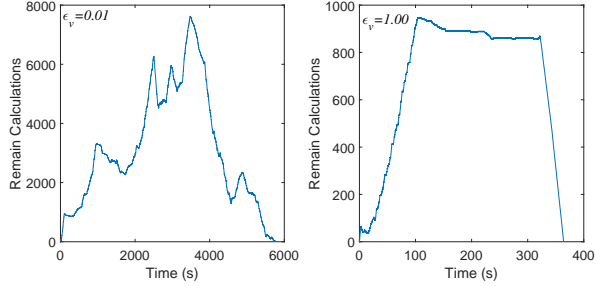}
	\caption{Number of hypercubes needed exploration during computation for 9-bus system with different values of $\varepsilon_V$}
	\label{fig-2}
\end{figure}

We present performance of the proposed method on the IEEE benchmark 9-bus system. The data can be found in Matpower toolbox \cite{Matpower}. By setting $\varepsilon_V = 10^{-2}$, $\varepsilon_R = 10^{-5}$, the ${\bf GPF_{sdp}}$ method finds 8 solutions in 5754.20 seconds. Results are shown in Table \ref{tab-4} and the number of hypercubes needed exploration during computation is illustrated in Figure \ref{fig-2}. It is noteworthy that if we set $\varepsilon_V = 1$, $\varepsilon_R = 10^{-5}$, the ${\bf GPF_{sdp}}$ method still finds all 8 solutions in 364 seconds. 
\begin{table}[ht]\renewcommand\arraystretch{1.0}
	\caption{PF Solutions of the 9-Bus System}
	\vspace{1pt}
	\centering
	\begin{tabular}{>{\scriptsize}c|>{\scriptsize}c|>{\scriptsize}c|>{\scriptsize}c|>{\scriptsize}c|>{\scriptsize}c|>{\scriptsize}c|>{\scriptsize}c|>{\scriptsize}c}
		\hline
		\multicolumn{9}{c}{\scriptsize Solution Index and Voltage Magnitude (p.u.)}\\
		\hline
		&	1 	&	2 	&	3 	&	4 	&	5 	&	6 	&	7 	&	8	\\
		\hline
		1	&	1.000 	&	1.000 	&	1.000 	&	1.000 	&   1.000 	&	1.000 	&	1.000 	&	1.000 	\\
		2	&	1.000 	&	1.000 	&	1.000 	&	1.000 	&	1.000 	&	1.000 	&	1.000 	&	1.000 	\\
		3	&	1.000 	&	1.000 	&	1.000 	&	1.000 	&	1.000 	&	1.000 	&	1.000 	&	1.000 	\\
		4	&	0.987 	&	0.857 	&	0.642 	&	0.653 	&	0.650 	&	0.593 	&	0.502 	&	0.468 	\\
		5	&	0.975 	&	0.767 	&	0.604 	&	0.706 	&	0.089 	&	0.134 	&	0.166 	&	0.139 	\\
		6	&	1.003 	&	0.681 	&	0.662 	&	0.901 	&	0.794 	&	0.585 	&	0.591 	&	0.785 	\\
		7	&	0.986 	&	0.077 	&	0.099 	&	0.810 	&	0.826 	&	0.108 	&	0.113 	&	0.756 	\\
		8	&	0.996 	&	0.582 	&	0.468 	&	0.779 	&	0.885 	&	0.539 	&	0.488 	&	0.775 	\\
		9	&	0.958 	&	0.709 	&	0.168 	&	0.121 	&	0.678 	&	0.487 	&	0.211 	&	0.195 	\\
		\hline
		\multicolumn{9}{c}{\scriptsize Solution Index and Voltage Angle (Degree)}\\
		\hline
		&	1 	&	2 	&	3 	&	4 	&	5 	&	6 	&	7 	&	8\\
		\hline
		1	& 0.000  & 0.000  & 0.000  & 0.000  & 0.000  & 0.000  & 0.000  & 0.000  \\
		2	& 9.669  & 13.121  & 38.188  & -5.605  & -11.276  & -3.312  & -10.853  & -66.510  \\
		3	& 4.771  & -8.555  & -16.594  & -11.017  & -21.901  & -74.141  & -78.898  & -76.575  \\
		4	& -2.407  & -6.785  & -10.479  & -9.257  & -11.126  & -13.039  & -15.204  & -16.338  \\
		5	& -4.017  & -13.577  & -21.705  & -16.850  & -70.677  & -82.206  & -82.693  & -92.533  \\
		6	& 1.926  & -12.752  & -20.908  & -14.184  & -25.498  & -79.026  & -83.730  & -80.215  \\
		7	& 0.622  & -61.709  & -51.702  & -17.085  & -24.208  & -94.276  & -102.230  & -80.621  \\
		8	& 3.799  & 3.040  & 25.603  & -13.121  & -17.883  & -14.210  & -22.902  & -74.063  \\
		9	& -4.350  & -11.705  & -46.153  & -62.001  & -21.313  & -27.181  & -56.573  & -83.433  \\
		\hline
	\end{tabular} 
	\label{tab-4} 
\end{table}

Next we test the ability of  ${\bf GPF_{sdp}}$ to detect infeasibility of PF problem. We investigate the following cases 
\begin{equation*}
\begin{lgathered}
\mbox{Case A} \left\{ \begin{lgathered}
V_k^{\min} = 0.90,\ V_k^{\max} = 1.10, \forall k \\
\Delta \theta _{ij}^{\max} = \pi / 9, \ \forall (ij)
\end{lgathered} \right.  \\
\mbox{Case B} \left\{ \begin{lgathered}
V_k^{\min} = 0.90,\ V_k^{\max} = 1.00, \forall k \\
\Delta \theta _{ij}^{\max} = \pi / 9, \ \forall (ij)
\end{lgathered} \right.  \\
\mbox{case C} \left\{ \begin{lgathered}
\begin{lgathered}
V_k^{\min} = 0.65, V_k^{\max} = 1.00 \\
\theta _k^{\min} = -\pi/12, \theta _k^{\max} = \pi / 12 \\
\end{lgathered}, \quad \forall k \not\in\{7,8 \} \\
\begin{lgathered}
V_k^{\min} = 0.00, V_k^{\max} = 1.00 \\
\theta _k^{\min} =-7\pi/18, \theta _k^{\max} = \pi / 12 \\
\end{lgathered}, \quad \forall k \in \{ 7,8 \} 
\end{lgathered} \right.  \\
\end{lgathered}
\end{equation*}

\begin{figure}
	\centering
	\includegraphics[scale=0.26]{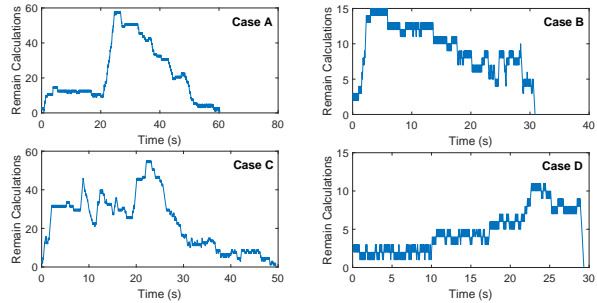}
	\caption{Number of hypercubes needed exploration during computation for Case A-Case D}
	\label{fig-3}
\end{figure}

In Case D, we intentionally modify the load data by multiplying the active power demand with a constant vector $\lambda$, such that the system status  is near to a bifurcation point. By setting $\varepsilon_V = 10^{-2}$, $\varepsilon_R = 10^{-5}$, the results are shown in Table \ref{tab-5}. The number of hypercubes needed exploration during computation is illustrated in Figure \ref{fig-3}. From the results of former 3 cases, we can see that when a smaller searching region is specified, the computation time can be reduced and feasibility does not have to be assumed priorly. As for Case D, when $\lambda > 2.52227$, the system has no power flow solution. Our method successfully detects infeasibility in 0.1 seconds. When $\lambda = 2.52226$, Algorithm 1 returns 2 solutions that are very close to each other with 29.33 seconds. The number of hypercubes needed exploration is illustrated in Figure \ref{fig-3}. This case indicates our method can efficiently detect infeasible power demand. 
\begin{table}[ht]\renewcommand\arraystretch{1.0}
	\caption{Results of detecting solutions of 9-bus system with restricted region}
	\vspace{6pt}
	\centering
	\begin{tabular}{>{\scriptsize}c|>{\scriptsize}c|>{\scriptsize}c}
		\hline
		& Computation time (s) &  Number of solutions \\
		\hline
		Case A &  60.29   &   1 (solution 1 in Table \ref{tab-4}) \\
		\hline
		Case B &  30.85   &   0   \\
		\hline
		Case C &  49.57   &   1 (solution 2 in Table \ref{tab-4}) \\
		\hline
	\end{tabular} 
	\label{tab-5} 
\end{table}

\subsection{14-Bus System}
For IEEE benchmark 14-bus system, the data of which can also be found in Matpower toolbox, we set $\varepsilon_V = 10^{-2}$, $\varepsilon_R = 10^{-5}$ and only study the ability of  ${\bf GPF_{sdp}}$ to detect infeasibility of PF problem. We here investigate the following cases 
\begin{equation*}
\begin{lgathered}
\mbox{Case E} \left\{ \begin{lgathered}
V_k^{\min} = 0.90,\ V_k^{\max} = 1.00, \forall k \\
\Delta \theta _{ij}^{\max} = \pi / 9, \ \forall (ij)
\end{lgathered} \right.  \\
\mbox{Case F} \left\{ \begin{lgathered}
V_k^{\min} = 1.00,\ V_k^{\max} = 1.10, \forall k \\
\Delta \theta _{ij}^{\max} = \pi / 9, \ \forall (ij)
\end{lgathered} \right.  \\
\end{lgathered}
\end{equation*}

For Case E, it takes 1252.64 seconds and the infeasibility is reported, i.e. no power flow solution is found. For case F, one power flow solution is found by solving ${\bf GPF_{sdp}}$ after 1166.27 seconds, which is the same as the one computed purely by traditional power flow algorithms in Matpower. These two cases again demonstrate our method can effectively declare the infeasibility or present the feasible power flow solution within a more restricted region.

\section{Conclusions}
A convex optimization based method is proposed to locate all real solutions or detect infeasibility for the power flow problem. It relies on well-developed SDP or SOCP solver and provides rigorous certification that either all solutions have been found, or there is no solution at all. The downside of the method is its scalability. We admit that at current stage, there is no method that is able to compute all power flow solutions of large-scale power systems. Despite this, our method is still capable of locating the areas that may contain a solution by using a relative large $\varepsilon _V $. Parallel computing technique and high-performance SDP or SOCP solvers may further improve the efficiency of the proposed algorithm.

%

\end{document}